\documentclass[reqno]{amsart}
\usepackage{amssymb,latexsym}
\usepackage{hyperref}
\usepackage{color}

\newtheorem{thm}{Theorem}[section]

\newtheorem{cor}{Corollary}[section]
\newtheorem{prop}{Proposition}[section]

\theoremstyle{definition}

\theoremstyle{remark}

\makeatletter
\@addtoreset{figure}{section}
\def\thefigure{\thesection.\@arabic\c@figure}
\def\fps@figure{h, t}
\@addtoreset{table}{bsection}
\def\thetable{\thesection.\@arabic\c@table}
\def\fps@table{h, t}
\@addtoreset{equation}{section}

\makeatother

\allowdisplaybreaks

\begin{document}

\title[Euler and second-grade fluid Lagrangian flow]
{Smooth global Lagrangian flow for the 2D Euler and second-grade 
fluid equations}

\author[S. Shkoller]{Steve Shkoller}
\address[Steve Shkoller]{Department of Mathematics\\
University of California \\
Davis, CA 95616}
\email{shkoller@math.ucdavis.edu}

\date{April 20, 2000; current version May 29, 2000}

\begin{abstract}
We present a very simple proof of the global existence of a $C^\infty$ 
Lagrangian flow map for the 2D
Euler and second-grade fluid equations (on a compact Riemannian manifold
with boundary) which has $C^\infty$ dependence on initial data $u_0$ in
the class of $H^s$ divergence-free vector fields for $s>2$.
\end{abstract}

\maketitle

\section{Incompressible Euler equations}Let $(M,g)$ be a $C^\infty$ compact
oriented Riemannian $2$-manifold with smooth boundary ${\partial M}$,  let
$\nabla$ denote the Levi-Civita covariant derivative, and let
$\mu$ denote the Riemannian volume form.  The incompressible
Euler equations are given by
\begin{equation}\label{E}
\begin{array}{c}
\partial_tu + \nabla_uu = \operatorname{grad}p\\
\operatorname{div}u=0, \ \ u(0)=u_0,  \ \ g(u,n)=0 \text{ on } {\partial M},
\end{array}
\end{equation}
where $p(t,x)$ is the pressure function, determined (modulo constants) by
solving the Neumann problem $-\triangle p = \operatorname{div}\nabla_uu$
with boundary condition $g(\operatorname{grad}p, n)= S_n(u)$, $S_n$ denoting
the second-fundamental form of ${\partial M}$.

The now standard global existence result for two-dimensional classical 
solutions 
states that for initial data $u_0 \in \chi^s \equiv \{ v \in H^s(TM) \ | \
\operatorname{div}u=0, \ g(u,n)=0\}$, $s>2$, the solution $u$ is in
$C^0({\mathbb R},\chi^s)$ and has $C^0$ dependence on $u_0$ (see, for example,
Taylor's book \cite{T}).
Equation (\ref{E}) gives the Eulerian or spatial representation of the dynamics
of the fluid.  The Lagrangian representation which is in terms of the
volume-preserving fluid
particle motion or flow map $\eta(t,x)$ is obtained by solving
\begin{equation}\label{flow}
\begin{array}{c}
\partial_t\eta(t,x) = u(t,\eta(t,x)),\\
\eta(0,x)=x.
\end{array}
\end{equation}
This is an ordinary differential equation on the infinite dimensional
volume-preserving diffeomorphism group ${\mathcal D}_\mu^s$, the set
of $H^s$ class bijective maps of $M$ into itself with $H^s$ inverses which
leave ${\partial M}$ invariant.
Ebin \& Marsden \cite{EM} proved that ${\mathcal D}_\mu^s$ is a 
$C^\infty$ manifold whenever $s>2$.  They also showed that for an
interval $I$, whenever $u\in C^0(I,\chi^s)$ and $s>3$, there exists 
a unique solution $\eta \in C^1(I,{\mathcal D}_\mu^s)$ to (\ref{flow}).
Thus, for $s>3$ the existence of a global $C^1$ flow map immediately follows
from the fact that $u$ remains bounded in $H^s$ for all time.
 It is often essential, however, for the Euler flow
to depend smoothly on the initial data;  in the case of vortex methods, for
example, Hald in Assumption 3 of \cite{H} requires this as a necessary 
condition to establish convergence.  

\begin{thm}\label{thm1}
For $u_0 \in \chi^s$, $s>2$, there exists a unique global solution to 
(\ref{L})
which is in $C^\infty({\mathbb R},T{\mathcal D}_\mu^s)$ and has $C^\infty$ 
dependence on $u_0$.
\end{thm}
\begin{proof}
The smoothness of the flow map follows by considering the Lagrangian version
of (\ref{E}) given  by
\begin{equation}\label{L}
\begin{array}{c}
\frac{D}{dt} \partial_t \eta(t,x) = -\operatorname{grad} p(t,\eta(t,x)),
 \ \ \ \operatorname{det} T\eta(t,x)=1,\\
\partial_t\eta(0,x) = u_0(x),\\
\eta(0,x)=x,
\end{array}
\end{equation}
where $T\eta(t,x)$ denotes the tangent map of $\eta$ (which in local coordinates
is given by the 2x2 matrix of partial derivatives 
$\partial \eta^i/\partial x^j$), and where $D/dt$ is the
covariant derivative along the curve $t \mapsto\eta(t,x)$ (which in 
Euclidean space is the usual partial time derivative).  Since
$$\operatorname{grad}p \circ \eta = 
\operatorname{grad}\triangle^{-1} \left[
\operatorname{Tr}(\nabla u \cdot \nabla u) +
\operatorname{Ric}(u,u)\right] \circ \eta,$$
where $\operatorname{Ric}$ is the Ricci curvature of $M$, and since
$S_n$ is $C^\infty$ and
$H^{s-1}(TM)$ forms a multiplicative algebra whenever $s>2$, we 
see that the linear operator 
$ u \mapsto \operatorname{grad}\triangle^{-1} [
\operatorname{Tr}(\nabla u \cdot \nabla u)$$+$$
\operatorname{Ric}(u,u)]$ maps $H^s$ back into $H^s$.  Denote by
$f:T{\mathcal D}_\mu^s \rightarrow TT{\mathcal D}_\mu^s$ the vector field
$$(\eta, \partial_t\eta) \mapsto
\operatorname{grad}\triangle^{-1} \left[
\operatorname{Tr}(\nabla u \cdot \nabla u) +
\operatorname{Ric}(u,u)\right] \circ \eta.$$
Then, 
$$f(\eta, \partial_t \eta) = 
\operatorname{grad}_\eta\triangle_\eta^{-1} \left[
\operatorname{Tr}(\nabla_\eta \partial_t\eta \cdot \nabla_\eta \partial_t\eta) 
+ \operatorname{Ric_\eta}(\partial_t\eta,\partial_t\eta)\right],$$
where $\operatorname{grad}_\eta g = [\operatorname{grad}(g \circ \eta^{-1})]
 \circ\eta$ for all $g\in H^s(M)$,
$\operatorname{div}_\eta X_\eta = [\operatorname{div}(X_\eta \circ 
\eta^{-1})] \circ\eta$ and $\nabla_\eta(X_\eta) = 
[\nabla (X_\eta \circ \eta^{-1}] \circ \eta$ for all $X_\eta \in 
T_\eta{\mathcal D}_\mu^s$,  $\triangle_\eta = \operatorname{div}_\eta
\circ \operatorname{grad}_\eta$, and $\operatorname{Ric_\eta}= 
\operatorname{Ric} \circ \eta$.  It follows from Lemmas 4,5, and 6 in
\cite{S} and Appendix A in \cite{EM} that $f$ is a $C^\infty$ vector field.
Thus (\ref{L}) is an ordinary 
differential equation on the tangent bundle $T{\mathcal D}_\mu^s$ governed
by a $C^\infty$ vector field on $T{\mathcal D}_\mu^s$; it immediately follows
from the fundamental theorem of ordinary differential equations on
Hilbert manifolds, that (\ref{L}) has a unique $C^\infty$ solution on 
{\it finite} time
intervals which depends smoothly on the initial velocity field $u_0$,
i.e., there exists a unique solution $\partial_t \eta
\in C^\infty( (-T,T), T{\mathcal D}_\mu^s)$ with $C^\infty$ dependence on
initial data $u_0$, where $T$ depends only on $\|u_0\|_{H^s}$.

When $s>3$, this interval can be extended globally to ${\mathbb R}$ by
virtue of $\eta$ remaining in ${\mathcal D}_\mu^s$.  Unfortunately, the
global existence and uniquess of a $C^\infty$ flow map $\eta(t,x)$ does
not follow for initial data $u_0 \in \chi^s$ for $s \in (2,3]$, so we
provide a simple argument to fill this gap.
We must show
that $\eta$ can be continued in ${\mathcal D}_\mu^s$.  It suffices to prove
that $T\eta$ and $T\eta^{-1}$ are both bounded in $H^{s-1}$.  This is easily
achieved using energy estimates.  We have that
$$\frac{D}{dt} T \eta = \nabla \partial_t\eta = \nabla u \cdot T\eta$$
and
$$\frac{D}{dt} T \eta^{-1} = -T\eta^{-1} \cdot \nabla \partial_t\eta \cdot
T\eta^{-1} = -T\eta^{-1} \cdot \nabla u.$$
Computing the $H^{s-1}$ norm of $T\eta$ and $T\eta^{-1}$, respectively, we
obtain
$$ {\frac{1}{2}}\frac{d}{dt} \| T\eta\|_{H^{s-1}}
= \langle D^{s-1}(\nabla u \cdot T\eta), D^{s-1} T\eta \rangle_{L^2},$$
and
$$ {\frac{1}{2}}\frac{d}{dt} \| T\eta^{-1}\|_{H^{s-1}}
= \langle D^{s-1}(T\eta^{-1} \cdot\nabla u),D^{s-1} T\eta^{-1}\rangle_{L^2}.$$
It is easy to estimate
\begin{align*}
\langle D^{s-1}(\nabla u \cdot T\eta), D^{s-1} T\eta \rangle_{L^2} & \le 
C( \|\nabla u\|_{L^\infty}  \|T\eta\|^2_{H^{s-1}} +
 \|\nabla u\|_{H^{s-1}}  \|T\eta\|_{L^\infty}   \|T\eta\|_{H^{s-1}}) \\
&\le C( \| \nabla u \|_{L^\infty}  \|T\eta\|^2_{H^{s-1}} +
 \|u\|_{H^{s}}  \|T\eta\|^2_{H^{s-1}}) 
\end{align*}
where the first inequality is due to Cauchy-Schwartz and  Moser's inequalities 
and the second is the Sobolev embedding theorem. Similarly,
$$
\langle D^{s-1}(-T \eta^{-1} \cdot \nabla u ), D^{s-1} T\eta^{-1} \rangle_{L^2} 
\le C( \| \nabla u \|_{L^\infty}  \|T\eta^{-1}\|^2_{H^{s-1}} +
 \|u\|_{H^{s}}  \|T\eta^{-1}\|^2_{H^{s-1}}).
$$
Since the solution $u$ to (\ref{E}) is in $\chi^s$ for all $t$, we have
that $\|u\|_{H^s}$ is bounded for all $t$.  Because the vorticity $\omega=
\operatorname{curl}u$ is in $L^\infty$, we have by Lemma 2.4 in Chapter 17 of
\cite{T} that $\| \nabla u\|_{L^\infty} \le C(1+\log\|u\|_{H^s})$; hence
$\| \nabla u\|_{L^\infty}$ is bounded for $t$.  It then follows that
$\eta$ and $\eta^{-1}$ are in ${\mathcal D}_\mu^s$ for all time.
\end{proof}

\section{Second-grade fluid equations} 
In this section, we establish the global existence of a $C^\infty$  Lagrangian
flow map for the second-grade fluids equations, also known as the
isotropic averaged Euler or Euler-$\alpha$ equations, which has $C^\infty$ 
dependence on intial data. 
 These equations are given on $(M,g)$ by
\begin{equation}\label{E2}
\begin{array}{c}
\partial_t(1-\alpha\triangle_r)u - \nu \triangle_r u +
\nabla_u(1-\alpha\triangle_r)u
-\alpha (\nabla u)^t \cdot \triangle_r u = -\text{\text{grad }}p ,\\
\operatorname{div} u =0, \ \ u(0)=u_0, \ \ u=0 \text{ on } {\partial M}\\
\alpha >0, \nu \ge 0 \ \ \ \triangle_r = -(d \delta + \delta d) + 2\text{Ric},
\end{array}
\end{equation}
(see \cite{S}),
and were first derived in 1955 by Rivlin\&Ericksenn \cite{RE} in Euclidean
space (Ric$=0$) as a first-order correction to  the Navier-Stokes equations.
In Euclidean space the operator $\triangle_r$ is just the component-wise
Laplacian, and the equation may be written as
$$\partial_t (1-\alpha \triangle) u - \nu\triangle u +
\text{curl}(1-\alpha\triangle)u \times u = -\text{grad } p.$$

For convenience, we set $\alpha =1$. We define the
unbounded, self-adjoint operator $(1-{\mathcal L})=(1-2\text{Def}^*\text{Def})$ 
on $L^2(TM)$ with domain $H^2(TM)\cap H^1_0(TM)$.  The operator $\operatorname
{Def}^*$ is the formal adjoint of $\operatorname{Def}$ with respect to $L^2$;
$2\text{Def}^* \text{Def }u = -(\triangle+\text{grad }\text{div}+ 2 \text{Ric})u$
so that
$2\text{Def}^* \text{Def }u = -(\triangle + 2 \text{Ric})u$ if $\text{div }u=0$.
We let ${\mathcal D}_{\mu,D}^s$ denote the subgroup of ${\mathcal D}_\mu^s$
whose elements restrict to the identity on the boundary ${\partial M}$.  
${\mathcal D}_{\mu,D}^s$ is a $C^\infty$ manifold (see \cite{EM} and \cite{S}).
Let $\chi^s_D = \{ u \in \chi^s \ | \ u=0 \text{ on } {\partial M}\}$.

The following is Proposition 5 in \cite{S}.
\begin{prop}\label{geodesic}
For $s>2$, let $\eta(t)$ be a curve in ${\mathcal D}_{\mu,D}^s$,
and set $u(t)=\partial_t\eta \circ
\eta(t)^{-1}$.  Then $u$ is a solution of the initial-boundary value problem
(\ref{E2}) with Dirichlet boundary conditions $u=0$ on ${\partial M}$
if and only if
\begin{equation}\label{L2}
\begin{array}{c}
\overline{\mathcal P}_\eta \circ \left[
\frac{\nabla \dot \eta}{dt} + \left[-\nu(1-{\mathcal L})^{-1}\triangle_r u +
{\mathcal U}(u) + {\mathcal R}(u)\right]
\circ \eta \right] =0, \ \ \operatorname{Det}T\eta(t,x)=1,\\
\partial_t\eta(0,x)=u_0(x),\\
\eta(0,x)=x,\\
\end{array}
\end{equation}
where
\begin{align*}
{\mathcal U}(u) =& (1-{\mathcal L})^{-1}\bigl\{
\text{\rm div}\left[ \nabla u \cdot
\nabla u^t + \nabla u \cdot \nabla u - \nabla u^t \cdot \nabla u\right]
+\text{\rm \text{grad }Tr}(\nabla u \cdot \nabla u)\bigr\} \\
{\mathcal R}(u) =& (1-{\mathcal L})^{-1} \bigr\{  \text{\rm Tr}
\left[ \nabla \left( R(u,\cdot)u \right) +R(u,\cdot) \nabla u +
R(\nabla u, \cdot)u \right] \\
&\qquad\qquad\qquad  + \text{\rm \text{grad }Ric}(u,u)
-(\nabla_u\text{\rm Ric}) \cdot u + \nabla u^t \cdot \text{\rm Ric}(u)
\bigr\}  ,
\end{align*}
and $\overline{\mathcal P}_\eta: T_\eta{\mathcal D}_D^s \rightarrow
T_\eta{\mathcal D}_{\mu,D}^s$ is the Stokes projector defined by
\begin{equation}\nonumber
\begin{array}{c}
\overline{\mathcal P}_\eta: T_\eta{\mathcal D}_{\mu,D}^s \rightarrow
T_\eta{\mathcal D}_{\mu,D}^s,\\
\overline{{\mathcal P}}_\eta(X_\eta) = \left[ {\mathcal P}_e(X_\eta\circ
\eta^{-1})\right] \circ\eta,
\end{array}
\end{equation}
and where ${\mathcal P}_e(F)=v$, $v$ being the unique solution of the Stokes 
problem
\begin{equation}\nonumber
\begin{array}{c}
(1-{\mathcal L}) v + \text{\rm grad }p = (1-{\mathcal L})F, \\
\text{\rm div }v =0,\\
v =0 \text{ \rm on } {\partial M}.
\end{array}
\end{equation}
\end{prop}
\noindent Equation (\ref{L2}) 
is an ordinary differential equation for the Lagrangian
flow.  Notice again that $H^{s-1}$, $s>2$, forms a 
multiplicative algebra, so that both ${\mathcal U}$ and ${\mathcal R}$ map
$H^s$ into $H^s$.

\begin{thm}\label{thm2}
For $u_0 \in \chi_D^s$, $s>2$,  and $\nu\ge 0$, there exists a unique global 
solution to (\ref{L2})
which is in $C^\infty({\mathbb R},T{\mathcal D}_\mu^s)$ and has $C^\infty$ 
dependence on $u_0$.
\end{thm}

We note that one cannot prove the statement of this theorem from an
analysis of (\ref{E2}) alone (see \cite{CG} and \cite{GGS}, and references
therein). 

\begin{proof}
The ordinary differential equation (\ref{L2}) can be written
as $\partial_{tt} \eta = S(\eta,\partial_t\eta)$ (see page 23 in \cite{S}).
Remarkably, $S:T{\mathcal D}_{\mu,D}^s \rightarrow TT{\mathcal D}_{\mu,D}^s$
is a $C^\infty$ vector field, and [\cite{S}, Theorem 2] provides
the existence of a unique short-time solution to (\ref{L2}) in $C^\infty((-T,T),
T{\mathcal D}_{\mu,D}^s)$ which depends smoothly on $u_0$, and where
$T$ only depends on $\|u_0\|_{H^s}$.

Thus, it suffices to prove that the solution curve $\eta$ does not
leave ${\mathcal D}_{\mu,D}^s$.  Following the proof of Theorem \ref{thm1},
and using the fact that the solution $u(t,x)$ to (\ref{E2}) remains in
$H^s$ for all time (\cite{CG,GGS}), it suffices to prove that $\nabla u$
is bounded in $L^\infty$.

Letting $q = \operatorname{curl}(1-\alpha\triangle_r) u$ denote the potential
vorticity, and computing the 
curl of (\ref{E2}), we obtain the 2D vorticity  form as
$$ \partial_t q + g(\text{grad }q,u) =\nu \text{curl }u.$$
It follows that for all $\nu \ge 0$, $q(t,x)$ is bounded in $L^2$ (conserved
when $\nu =0$) and therefore by standard elliptic
estimates $\nabla u(t,x)$ is bounded in $H^2$, and hence in $L^\infty$.
\end{proof}

As a consequence of Theorem \ref{thm2} being independent of viscosity,
we immediately obtain the following:
\begin{cor}
Let $\eta^\nu(t,x)$ denote the Lagrangian flow solving (\ref{L2}) for 
$\nu >0$, so that
$u^\nu=\partial_t \eta^\nu \circ {\eta^\nu}^{-1}$ solves (\ref{E2}).  Then for
$u_0 \in \chi_D^s$,  $s>2$, the viscous solution $\eta^\nu \in 
C^\infty({\mathbb R}, T{\mathcal D}_\mu^s)$ converges regularly (in $H^s$) 
to the inviscid solution $\eta^0\in C^\infty({\mathbb R}, T{\mathcal D}_\mu^s)$.
Consequently $u^\nu \rightarrow u^0$ in $H^s$ on infinite-time intervals.
\end{cor}

This gives an improvement of  Busuioc's result in \cite{B} in two ways: 1) we
are able to prove the regular limit of zero viscosity on manifolds with
boundary, and 2) in the Lagrangian framework, we are able to get $C^\infty$
in time solutions.

\section*{Acknowledgments}
Research  was partially supported by the NSF-KDI grant ATM-98-73133
and the Alfred P. Sloan Foundation Research Fellowship.

\end{document}